\documentclass[12pt]{article}
\usepackage{amsmath}
\usepackage{algorithm}
\usepackage{algorithmicx}
\usepackage{algpseudocode}
\usepackage{lineno}
\usepackage{bm}
\usepackage{t1enc}
\usepackage[latin1]{inputenc}
\usepackage[english]{babel}
\usepackage{amsmath,amsthm}
\usepackage{amsfonts}
\usepackage{latexsym}
\usepackage[dvips]{graphicx}
\usepackage{graphicx}
\usepackage[natural]{xcolor}
\usepackage{enumerate}
\usepackage{tikz}
\usetikzlibrary{shapes.geometric, arrows}
\usepackage{indentfirst}
\usepackage{listings}
\usepackage{booktabs}
\usepackage{natbib,stfloats}
\usepackage{mathrsfs}
\usepackage{amssymb}
\usepackage{authblk}
\allowdisplaybreaks 
\newtheorem{theorem}{Theorem}[section]

\usepackage{url}
\allowdisplaybreaks 

\begin{document}

\title{PASE: A Massively Parallel Augmented Subspace Eigensolver for Large Scale Eigenvalue Problems\footnote{This work was supported by the National Key Research and Development 
Program of China (2019YFA0709601, 2023YFB3309104), 
National Natural Science Foundations of China (NSFC 1233000214), Science Challenge Project (TZ2024009), 
Beijing Natural Science Foundation (Z200003), 
the National Center for Mathematics and Interdisciplinary Science, CAS.}}



\author[1]{Yangfei Liao\thanks{Email: \tt liaoyangfei@lsec.cc.ac.cn}}
\author[2]{Haochen Liu\thanks{Email: \tt liuhaochen@lsec.cc.ac.cn}}
\author[3]{Hehu Xie\thanks{Email: \tt hhxie@lsec.cc.ac.cn}}
\author[4]{Zijing Wang$^{*}$\thanks{Email: \tt zjwang@lsec.cc.ac.cn} \thanks{$^{*}$Corresponding author.}}

\affil[1234]{SKLMS, ICMSEC, Academy of Mathematics and Systems Science, Chinese Academy of Sciences, Beijing 100190, China, and School of Mathematical Sciences, University of Chinese Academy of Sciences, Beijing 100049, China}


\date{}


\maketitle


\abstract{
In this paper, we present a novel parallel augmented subspace method 
and build a package Parallel Augmented Subspace Eigensolver (PASE) for
solving large scale eigenvalue problems by the massively parallel finite element discretization.
Based on the augmented subspace, solving high dimensional eigenvalue problems can be
transformed to solving the corresponding linear equations and low dimensional eigenvalue problems on the augmented subspace.  
Thus the complexity of solving the eigenvalue problems by augmented subspace method will be comparable to that of solving the same dimensinal linear equations.
In order to improve the scalability and efficiency, we also present some implementing techniques for the
parallel augmented subspace method.  Based on parallel augmented subspace method and
the concerned implementing techniques, a package PASE is built for solving large scale eigenvalue problems.
Some numerical examples are provided to validate the efficiency and scalability of the proposed numerical methods.
}

\maketitle

\section{Introduction}
\label{sec_intro}
Eigenvalue problems arise in a large number of disciplines of science and engineering.
They act as the basic tool for build designing, bridges stability analysis,
structure dynamics, quantum chemistry, electrical networks, 
Markov chain techniques, chemical reactions,
data mining, big data, artificial intelligence, material sciences and so on.
Then the numerical computation of eigenpairs of large matrices is a problem
of major importance in many scientific and engineering applications.
Along with the development of modern science and engineering, there appear more and
more large scale eigenvalue problems. This brings the strong demand for efficient solver
for large scale eigenvalue problems \cite{Bai,Saad1,SLEPC}. Among these eigenvalue problems,
the ones from the discretization of differential operators act 
an important role in physics and engineering.
The appearance of high performance computers brings more ability for computing plenty 
of eigenpairs of
large scale eigenvalue problems \cite{SLEPC,LiWangXie, LiXieXuYouZhang,ZhangLiXieXuYou}.

Multigrid (MG) method has become one of popular and important solvers in
scientific computing due to its robustness and efficiency for solving large scale
unstructured linear systems of equations,
particularly when the matrix results from the discretization of some type of second-order
elliptic partial differential equations (PDE).
The key to the efficiency of multigrid method whether geometric or algebraic, is that the error
which can not be reduced by the smoothing (relaxation) process, can be eliminated by
the coarse grid correction process.  For more details and information, please refer to
\cite{BrambleZhang,BrennerScott,Hackbusch,Hackbusch_Book,McCormick,Shaidurov,Xu}.
The extension of the multigrid idea to the eigenvalue problems is natural in the field of
numerical algebra. Always, the multigrid method is adopted as the linear solver in solving
eigenvalue problems. This application can reduce the time cost to solve eigenvalue problems.
But this way can not improve the convergence rate of the eigensolver itself.
The essential reason for this difficulty is that even the linear eigenvalue problem is
actually a nonlinear equation.

Recently, a type of augmented subspace method is proposed for solving eigenvalue problems.
This method is designed based on an augmented subspace consisting of a low dimensional 
finite element space and current eigenfunction approximation in high dimensional finite space.
With the help of this special augmented subspace, we can transform solving
the large scale eigenvalue problem defined in the high dimensional finite element space
to solving the same scale linear problems and the small scale eigenvalue problem
defined in the low dimensional augmented subspace.
Actually, the augmented subspace method deduces a new
eigensolver which has the uniform convergence rate. And the convergence rate
depends on the approximation accuracy of the low dimensional subspace 
for the associated differential operator,
which can be improved easily by enlarging the low dimensional space 
(i.e., refine the mesh size of the coarse mesh).
The idea of the augmented subspace method is different from the methodology 
based on the Krylov subspace method.
The essential difference of the coarse finite element space from the Krylov 
subspace is that the former has the
approximation accuracy to the concerned differential operator in the operator-norm. 
The low dimensional subspace can be chosen as the finite element space defined on a coarse mesh
in the sequence of nested meshes.  
This understanding can bring more applications for solving eigenvalue problems 
which will also be investigated in this paper.  
Based on this consideration, we can design a type of
multigrid (or multilevel) method to solve the eigenvalue problems.
For more details, please refer to our previous papers \cite{ChenXieXu,DangWangXieZhou,
HongXieXu,Xie_JCP, LinXie_MultiLevel,Xie_IMA,XieZhangOwhadi, XuXieZhang,ZhangHanHeXieYou}.
Of course, the coarse subspace in the algebraic multigrid method can also act as the low dimensional subspace in the augmented subspace method.
For more details of this type of algebraic multigrid method for eigenvalue problems, please refer to
\cite{ZhangHanHeXieYou}.

The augmented subspace method makes the main computational work is 
to solve the linear equation in the high dimensional space.
It is well known that there exists many parallel numerical methods 
for solving the linear equations. With the help of augmented subspace, 
\cite{XuXieZhang} designs a type of eigenwise
parallel multilevel correction method for solving plenty of eigenpair approximations 
on the sequential distribution mesh.
The key for this type of parallel method is that, for different eigenvalues,
the corresponding boundary value problem and low-dimensional eigenvalue
problem can be solved in the parallel way since they are independent of each other
and there exists no data exchanging. Furthermore, this property means that we do
not need to do the inner products in the high dimensional spaces.
It is well known that decreasing inner product computation can improve the scalability
of the concerned numerical method.
Different from the previous work, the aim of this paper is to design a 
thoroughgoing parallel method for 
solving the eigenvalue problems based on the parallel distribution of mesh, 
augmented subspace method and parallel linear solver for the linear equations.
The bottleneck of the scalability for the parallel augmented subspace method is solving the
eigenvalue problem on the augmented subspace.
In order to enhance the scalability to a greater extent, we design some preprocessing techniques 
for the eigenvalue problem on the augmented subspace. This eigenvalue problem 
is solved by the General Conjugate Gradient Eigensolver (GCGE) which is a 
parallel computing package for solving eigenvalue problems 
\cite{LiWangXie,LiXieXuYouZhang,ZhangLiXieXuYou}.
We choose GCGE since it is built with the way of matrix-free and vector-free and suitable for
the special matrix and vector structure of the eigenvalue problem on the augmented subspace in this paper.
The second aim of this paper is to combine
this type of eigenpair parallel method to improve the scalability further for solving
plenty of eigenpair approximations on the parallel distribution of mesh.

An outline of the paper goes as follows. In Section 2, 
we introduce the augmented subspace method for the eigenvalue problems.
The parallel implementation of PASE will be designed in Section 3. In Section 4, we discuss how to address the scalability bottleneck in eigenvalue solvers using preconditioning techniques. Subsequently, we present the design of a batch scheme of PASE for solving a large number of eigenpairs in Section 5. 
In Section 6, some numerical examples
are provided to validate the efficiency, stability and scalability of PASE. 

\section{Augmented subspace method for eigenvalue problem}\label{sec_2}

In this section, we first provide a brief description of the finite element method for the second order elliptic eigenvalue problem. 
More importantly, we will elaborate the augmented subspace method which is an acceleration algorithm of the finite element method for eigenvalue problem. Additionally, the relevant error estimates of this method will be presented.

In this paper, we shall use the standard notation for Sobolev spaces $W^{s,p}(\Omega)$ and their associated norms and semi-norms (cf. \cite{Adams}). For $p=2$, we denote
$H^s(\Omega)=W^{s,2}(\Omega)$ and
$H_0^1(\Omega)=\{v\in H^1(\Omega):\ v|_{\partial\Omega}=0\}$,
where $v|_{\Omega}=0$ is in the sense of trace,
$\|\cdot\|_{s,\Omega}=\|\cdot\|_{s,2,\Omega}$.
In some places, $\|\cdot\|_{s,2,\Omega}$ should be viewed as piecewise
defined if it is necessary.

\subsection{The finite element method for eigenvalue problem}
For simplicity, we consider the following second order elliptic eigenvalue problem 
to illustrate the main idea of finite element method:
Find $(\lambda, u)$ such that
\begin{eqnarray}\label{LaplaceEigenProblem}
\left\{
\begin{array}{rcl}
-\nabla\cdot(\mathcal{A}\nabla u) &=&\lambda u, \quad {\rm in} \  \Omega,\\
u&=&0, \ \ \quad {\rm on}\  \partial\Omega,
\end{array}
\right.
\end{eqnarray}
where $\mathcal{A}$ is a symmetric and positive definite matrix with suitable
regularity, $\Omega\subset\mathcal{R}^d\ (d=2,3)$ is a bounded domain with
Lipschitz boundary $\partial\Omega$. We would like to point out the method in this
paper can be used to more general elliptic type of eigenvalue problems.

In order to use the finite element method to solve
the eigenvalue problem \eqref{LaplaceEigenProblem}, we need to define
the corresponding variational form as follows:
Find $(\lambda, u )\in \mathcal{R}\times V$ such that $a(u,u)=1$ and
\begin{eqnarray}\label{weak_eigenvalue_problem}
a(u,v)=\lambda b(u,v),\quad \forall v\in V,
\end{eqnarray}
where $V:=H_0^1(\Omega)$ and
\begin{eqnarray*}
a(u,v)=\int_{\Omega}\mathcal{A}\nabla u\cdot\nabla v d\Omega,
 \ \ \ \  \ \ b(u,v) = \int_{\Omega}uv d\Omega.
\end{eqnarray*}
The norms $\|\cdot\|_a$ and $\|\cdot\|_b$ are defined by
\begin{eqnarray*}
\|v\|_a=\sqrt{a(v,v)}\ \ \ \ \ {\rm and}\ \ \ \ \ \|v\|_b=\sqrt{b(v,v)}.
\end{eqnarray*}
It is well known that the eigenvalue problem \eqref{weak_eigenvalue_problem}
has an eigenvalue sequence $\{\lambda_j \}$ (cf. \cite{BabuskaOsborn_1989,Chatelin}):
$$0<\lambda_1\leq \lambda_2\leq\cdots\leq\lambda_k\leq\cdots,\ \ \
\lim_{k\rightarrow\infty}\lambda_k=\infty,$$ and associated
eigenfunctions
$$u_1, u_2, \cdots, u_k, \cdots,$$
where $a(u_i,u_j)=\delta_{ij}$ ($\delta_{ij}$ denotes the Kronecker function).
In the sequence $\{\lambda_j\}$, the $\lambda_j$ are repeated according to their
geometric multiplicity.

Now, let us define the finite element approximations of the problem
(\ref{weak_eigenvalue_problem}). First we generate a shape-regular
generated  mesh $\mathcal{T}_h$ of the computing domain $\Omega\subset \mathcal{R}^d\
(d=2,3)$ into triangles or rectangles for $d=2$ (tetrahedrons or
hexahedrons for $d=3$). The diameter of a cell $K\in\mathcal{T}_h$
is denoted by $h_K$ and the mesh size $h$ describes the maximal diameter of all cells
$K\in\mathcal{T}_h$. Based on the mesh $\mathcal{T}_h$, we can construct
a finite element space denoted by $V_h \subset V$. For simplicity, 
we set $V_h$ as the standard Lagrange type of finite element space which 
is defined as follows
\begin{eqnarray*}
V_h = \left\{ v_h \in C(\Omega)\ \big|\ v_h|_{K} \in \mathcal{P}_k,
\ \ \forall K \in \mathcal{T}_h\right\}\cap H^1_0(\Omega),
\end{eqnarray*}
where $\mathcal{P}_k$ denotes the polynomial set of degree no more than $k$. 

Based on the space $V_h$, the standard finite element scheme for eigenvalue
problem \eqref{weak_eigenvalue_problem} is:
Find $(\bar{\lambda}_h, \bar{u}_h)\in \mathcal{R}\times V_h$
such that $a(\bar{u}_h,\bar{u}_h)=1$ and
\begin{eqnarray}\label{Weak_Eigenvalue_Discrete}
a(\bar{u}_h,v_h)=\bar{\lambda}_h b(\bar{u}_h,v_h),\quad\ \  \ \forall v_h\in V_h.
\end{eqnarray}
From \cite{BabuskaOsborn_1989,BabuskaOsborn_Book}, the  discrete eigenvalue
problem \eqref{Weak_Eigenvalue_Discrete} has eigenvalues:
$$0<\bar{\lambda}_{1,h}\leq \bar{\lambda}_{2,h}\leq\cdots\leq \bar{\lambda}_{k,h}
\leq\cdots\leq \bar{\lambda}_{N_h,h},$$
and corresponding eigenfunctions
\begin{eqnarray*}
\bar{u}_{1,h}, \bar{u}_{2,h}, \cdots, \bar{u}_{k,h}, \cdots, \bar{u}_{N_h,h},
\end{eqnarray*}
where $a(\bar{u}_{i,h},\bar{u}_{j,h})=\delta_{ij}$, $1\leq i,j\leq N_h$ ($N_h$ is
the dimension of the finite element space $V_h$).
The error estimates for the first $k$ eigenpair approximations computed by finite element method can be found in \cite{BabuskaOsborn_Book}.

\subsection{The augmented subspace mehtod}
Next, we propose the augmented subspace method for eigenvalue problems. The main idea of the augmented subspace method is to transform solving
  the eigenvalue problem in the high dimensional finite element space into the solution
  of the corresponding linear boundary value problems in the same high dimensional
  finite element space and eigenvalue problems on a very low dimensional augmented subspace.
  Since solving eigenvalue problems appears only in the low dimensional augmented subspace and
  the main work is to solve the linear boundary value problem on the high dimensional parallel finite element
  space, the augmented subspace method can make solving eigenvalue problem
  be not significantly more expensive than solving the corresponding linear boundary value problems.
  
  In order to design the augmented subspace method, we first generate coarse mesh $\mathcal{T}_H$ with the mesh size $H$ and fine mesh $\mathcal{T}_h$
  with the mesh size $h$. $V_H$ and $V_h$ denote the concerned finite element space defined on $\mathcal T_H$ and $\mathcal T_h$. For simplicity, we assume the coarse space $V_H$ is a subspace of the high dimensional finite element space $V_h$.
  
  We first present the augmented subspace method for solving $k$ eigenpairs of (\ref{Weak_Eigenvalue_Discrete}).
  This method contains solving the auxiliary linear boundary value problems
  in $V_h$ and the eigenvalue problem on the augmented subspace $V_{H,h}$ which is built by $V_H$ and finite element functions in $V_h$.
For the positive integer $\ell$ which denotes the current iteration step and the given $k$ eigenfunction approximations $u_{1,h}^{(\ell)}$, $\cdots$, $u_{k,h}^{(\ell)}$, which are the approximations for the $k$ eigenfunctions $\bar u_{1,h}$, $\cdots$, $\bar u_{k,h}$ of (\ref{Weak_Eigenvalue_Discrete}), we can perform the following augmented subspace iteration step, as defined by Algorithm \ref{Algorithm_k}, to improve the accuracy of $u_{1,h}^{(\ell)}$, $\cdots$, $u_{k,h}^{(\ell)}$.

\begin{algorithm}[htbp!]
\caption{One correction step for $k$ eigenpairs}
\label{Algorithm_k}
\begin{algorithmic}[0]
\State Solve the following linear boundary value problems:
Find $\widehat{u}_{i,h}^{(\ell+1)}\in V_h$ such that
\begin{equation}\label{Linear_Equation_k}
a(\widehat{u}_{i,h}^{(\ell+1)},v_h) = \lambda_{i,h}^{(\ell)}b(u_{i,h}^{(\ell)},v_h),
\ \  \forall v_h\in V_h,\ \ \ i=1, \cdots, k.
\end{equation}

\State Define the augmented subspace $V_{H,h}^{(\ell+1)} = V_H + {\rm span}\{\widehat{u}_{1,h}^{(\ell+1)}, \cdots, \widehat u_{k,h}^{(\ell+1)}\}$ 
and solve the following eigenvalue problem:
Find $(\lambda_{i,h}^{(\ell+1)},u_{i,h}^{(\ell+1)})\in \mathcal{R}\times V_{H,h}^{(\ell+1)}$ such that $a(u_{i,h}^{(\ell+1)},u_{i,h}^{(\ell+1)})=1$ and
\begin{equation}\label{Aug_Eigenvalue_Problem_k}
a(u_{i,h}^{(\ell+1)},v_{H,h}) = \lambda_{i,h}^{(\ell+1)}b(u_{i,h}^{(\ell+1)},v_{H,h}),
\ \  \forall v_{H,h}\in V_{H,h}^{(\ell+1)},\ \  i=1, \cdots, k.
\end{equation}\\
Use the notation below to represent the total effect of the above two steps:
\begin{equation*}
    (\lambda_{i,h}^{(\ell+1)}, u_{i,h}^{(\ell+1)}) = \text{Correction}(V_H, V_h, \lambda_{i,h}^{(\ell)}, u_{i,h}^{(\ell)}).
\end{equation*}
\end{algorithmic}
\end{algorithm}

From \cite{DangWangXieZhou}, the augmented subspace method defined by Algorithm \ref{Algorithm_k} has the following error estimates. The notation remains the same as \cite{DangWangXieZhou}. 

\begin{theorem}(\cite{DangWangXieZhou})\label{Theorem_Error_Estimate_k}
Let us define the spectral projection $F_{k,h}^{(\ell)}: V\mapsto {\rm span}\{u_{1,h}^{(\ell)}, \cdots, u_{k,h}^{(\ell)}\}$ for any integer $\ell \geq 1$ as follows:
\begin{eqnarray*}
a(F_{k,h}^{(\ell)}w, u_{i,h}^{(\ell)}) = a(w, u_{i,h}^{(\ell)}), \ \ \ i=1, \cdots, k\ \ {\rm for}\ w\in V.
\end{eqnarray*}
Then the exact eigenfunctions $\bar u_{1,h},\cdots, \bar u_{k,h}$ of (\ref{Weak_Eigenvalue_Discrete}) and the eigenfunction approximations $u_{1,h}^{(\ell+1)}$, $\cdots$,  $u_{k,h}^{(\ell+1)}$ from Algorithm \ref{Algorithm_k} with the integer $\ell > 1$ have the following error estimate:
\begin{eqnarray*}\label{Error_Estimate_Inverse}
 \left\|\bar u_{i,h} - F_{k,h}^{(\ell+1)}\bar u_{i,h} \right\|_a \leq
 \bar\lambda_{i,h} \sqrt{1+\frac{\eta_a^2(V_H)}{\bar\lambda_{1,h}\big(\delta_{k,i,h}^{(\ell+1)}\big)^2}}
\left(1+\frac{\bar\mu_{1,h}}{\delta_{k,i,h}^{(\ell)}}\right)\eta_a^2(V_H)\left\|\bar u_{i,h} - F_{k,h}^{(\ell)}\bar u_{i,h} \right\|_a,
\end{eqnarray*}
where $\delta_{k,i,h}^{(\ell)} $ is defined as follows:
\begin{eqnarray*}
\delta_{k,i,h}^{(\ell)} = \min_{j\not\in \{1, \cdots, k\}}\left|\frac{1}{\lambda_{j,h}^{(\ell)}}-\frac{1}{\bar\lambda_{i,h}}\right|,\ \ \ i=1, \cdots, k.
\end{eqnarray*}
Furthermore, the following $\left\|\cdot\right\|_b$-norm error estimate holds:
\begin{eqnarray*}
\left\|\bar u_{i,h} -F_{k,h}^{(\ell+1)}\bar u_{i,h} \right\|_b\leq 
\left(1+\frac{\bar\mu_{1,h}}{\delta_{k,i,h}^{(\ell+1)}}\right)\eta_a(V_H) \left\|\bar u_{i,h} -F_{k,h}^{(\ell+1)}\bar u_{i,h}\right\|_a.
\end{eqnarray*}
\end{theorem}

\begin{algorithm}
  \caption{Augmented Subspace Method for $k$ Eigenpairs}
  \label{Algorithm_PASE}
  \begin{algorithmic}[1]
    \State Solve the following eigenvalue problem: Find $(\lambda_{i,H}, u_{i,H}) \in \mathbb{R} \times V_H$ such that $a(u_{i,H}, u_{i,H}) = 1$ for $i = 1, \dots, k$ and $a(u_{i,H}, u_{j,H}) = 0$ for $i \neq j$:
    \begin{equation*}\label{initial_eigenproblem}
      a(u_{i,H}, v_H) = \lambda_{i,H} b(u_{i,H}, v_H), \quad \forall v_H \in V_H, \quad i = 1, \dots, k.
    \end{equation*}

    \State Define $\lambda^{(1)}_{i,h} = \lambda_{i,H}$, $u^{(1)}_{i,h} = u_{i,H}$ for $i = 1, \dots, k$ and set $\ell = 1$.

    \For{$\ell = 1,2,\dots$}
      \State Compute:
      \begin{equation*}
        (\lambda_{i,h}^{(\ell+1)}, u_{i,h}^{(\ell+1)}) = \text{Correction}(V_H, V_h, \lambda_{i,h}^{(\ell)}, u_{i,h}^{(\ell)}).
      \end{equation*}
      \If{the iteration step converges}
        \State Exit loop and return the result.
      \EndIf
    \EndFor

  \end{algorithmic}
\end{algorithm}

We use Algorithm \ref{Algorithm_PASE} to illustrate the complete process of the augmented subspace method, as it enables us to solve larger-scale problems with a greater number of eigenpairs.

\section{Parallel augmented subspace eigensolver}\label{cha_para}
In this section, we introduce the package PASE (parallel augmented subspace eigensolver), which is designed based on the parallel implementation of the augmented subspace method defined by Algorithm \ref{Algorithm_PASE}. 

As we mentioned in Section \ref{sec_intro}, the eigenwise parallel method in \cite{XuXieZhang} is based on the sequential distribution mesh. Although this method completely avoids the data exchange between processes, both the storage and computation of matrices were inherently sequential. Consequently, this approach is easily constrained by the memory of a single node and the number of process available. This limitation prevents any efficiency gains through parallelization of matrix computing.

Innovatively, the method proposed in this paper achieves a thoroughgoing parallelization, including parallel distribution mesh and parallel matrix operation. For Algorithm \ref{Algorithm_PASE}, parallel linear solver and eigensolver are employed to solve \eqref{Linear_Equation_k} and \eqref{Aug_Eigenvalue_Problem_k}, respectively. This approach enables us to solve larger-scale problems with a greater number of eigenpairs.

We use Algorithm \ref{Algorithm_PASE_Algebric} to demonstrate the algebraic version of algorithm which is actually used in PASE along with the detailed implementation, where the input matrices can be obtained from geometric or algebraic multigrid methods. It should be noted that the input interpolation operator in Algorithm \ref{Algorithm_PASE_Algebric} can involve multiplying by a single interpolation matrix or by a sequence of interpolation matrices, which depends on the choice of coarse and fine meshes in the multigrid structure.

\begin{algorithm}[htbp]
  \caption{Algebraic parallel augmented subspace method in PASE}
  \label{Algorithm_PASE_Algebric}
  \begin{algorithmic}[1]
  \Require Stiffness and mass matrices $A_H$, $B_H$ on the coarse finite element space, $A_h$, $B_h$ on the fine finite element space, interpolation operator $I_H^h$ and its transpose $I_h^H$.
  \Ensure The first $k$ eigenpairs $(\Lambda, \mathbf{u}_{h})$ such that $A_h\mathbf{u}_h=B_h\mathbf{u}_h\Lambda$, where $\Lambda = \text{diag}(\lambda_1, ..., \lambda_k)$ and $\mathbf{u}_h = [u_{1,h},...,u_{k,h}]$.
  
  \State Solve the initial eigenvalue problem:
  \begin{equation*}
  A_H \mathbf{u}_H = B_H \mathbf{u}_H\Lambda_H.
  \end{equation*}
  
  \State Interpolate $\mathbf{u}_H$ to obtain $\mathbf{u}_h^{(\ell)} = I_H^h\mathbf{u}_H$ and set $\Lambda_h^{(\ell)}=\Lambda_H$, $\ell=1$.

  \State Solve the linear equation with initial solution $\mathbf{u}_h^{(\ell)}$:
  \begin{equation}\label{presoomthing}
  A_h \widehat{\mathbf{u}}_h^{(\ell+1)} = B_h \mathbf{u}_h^{(\ell)}\Lambda_h^{(\ell)}.
  \end{equation}

  \State Assemble the composite matrices $A_{Hh}$ and $B_{Hh}$ defined by
  \begin{equation}\label{Definition_Hh}
    A_{Hh} =
    \begin{bmatrix}
    A_H & \mathbf{a}_h\\
    \mathbf{a}_h^\top  & \boldsymbol{\alpha}
    \end{bmatrix}, \quad
   B_{Hh} =
    \begin{bmatrix}
    B_H & \mathbf{b}_h\\
    \mathbf{b}_h^\top  & \boldsymbol{\beta}
    \end{bmatrix}.
  \end{equation}
  where $\mathbf{a}_h = I_h^HA_h\widehat{\mathbf{u}}_h^{(\ell+1)}$, $\boldsymbol{\alpha} = (\widehat{\mathbf{u}}_h^{(\ell+1)})^T A_h \widehat{\mathbf{u}}_h^{(\ell+1)}$, $\mathbf{b}_h=I_h^HB_h\widehat{\mathbf{u}}_h^{(\ell+1)}$, and $\boldsymbol{\beta}=(\widehat{\mathbf{u}}_h^{(\ell+1)})^T B_h \widehat{\mathbf{u}}_h^{(\ell+1)}$.

  \State Solve the following eigenvalue problem with initial guess $\mathbf{u}_H = \mathbf{0}$ and $\gamma = \mathbf{I}$:
  \begin{equation}\label{Augmented_GCGE}
  A_{Hh} \mathbf{u}_{Hh} = B_{Hh} \mathbf{u}_{Hh}\Lambda_h, 
  \quad \mathbf{u}_{Hh}\triangleq \left(
    \begin{array}{c}
    \mathbf{u}_H\\
   \boldsymbol{\gamma}
    \end{array}
    \right).
  \end{equation}

  \State Solve (\ref{Augmented_GCGE}) to obtain $\Lambda_h^{(\ell+1)}$,  
  $\mathbf{u}_{Hh}^{(l+1)}=\left(
    \begin{array}{c}
    \mathbf{u}_H^{(l+1)}\\
   \boldsymbol{\gamma}^{(\ell+1)}
    \end{array}
    \right)$ and compute $\widehat{\mathbf{u}}_h^{(\ell+1)}\leftarrow I_H^h\mathbf{u}_H^{(\ell+1)} +\boldsymbol{\gamma} \widehat{\mathbf{u}}_h^{(\ell+1)}$.

  \State Solve the following linear equation with the initial guess $\widehat{\mathbf{u}}_h^{(\ell+1)}$ to obtain $\mathbf{u}_h^{(\ell+1)}$:
  \begin{equation}\label{postsoomthing}
  A_h \mathbf{u}_h^{(\ell+1)} = B_h \widehat{\mathbf{u}}_h^{(\ell+1)}\Lambda_h^{(\ell+1)}.
  \end{equation}

  \State Set $\ell=\ell+1$ and go to Step 3 for next iteration until convergence.

  \end{algorithmic}
\end{algorithm}

In Algorithm \ref{Algorithm_PASE_Algebric}, the matrices $A_H$, $B_H$, $A_h$, $B_h$, interpolation operator $I_H^h$ and restriction operator $I_h^H$ are all sparse and parallelized, while the corresponding eigenvectors $u_H$ and $u_h$ are also stored as parallel multi-vectors. Consequently, We have thoroughly parallelized most of the computing steps such as Step 1 for computing initial eigenpairs, Step 3 for solving linear systems with iterative methods, convergence check, interpolation and so on.

Next, we will focus more on discussing the eigenvalue problem \eqref{Aug_Eigenvalue_Problem_k} and its matrix version \eqref{Augmented_GCGE} on the augmented space, which is one of the most important parts of the algorithm. Let $N_H$ and $\{\phi_{i,H}\}_{1\leq i\leq N_H}$ denote the dimension and the Lagrange basis functions for the coarse finite element space $V_H$. Suppose we already have an approximation of the $k$ eigenfunctions on $V_h$ denoted by $\{\widehat{u}_{j,h}\}_{1\leq j\leq k}$, on the fine space. Then the function in $V_{H,h}=V_H + \text{span}\{\widehat{u}_{1,h},...,\widehat{u}_{k,h}\}$ can be expressed as $u_{H,h}=\sum_{i=1}^{N_H}u_i\phi_{i,H}+\sum_{j=1}^k \gamma_j\widehat{u}_{j,h}$. And we can also derive the composite part of the algebric eigenvalue problem \eqref{Augmented_GCGE} as:

\begin{equation*}
\begin{split}
&\mathbf{a}_h = \big[a(\phi_{i,H},\widehat{u}_{j,h})\big]_{\scriptscriptstyle 1\leq i\leq N_H, 1\leq j\leq k} \in\mathcal{R}^{N_H\times k}, \\
&\boldsymbol{\alpha} =\big[a(\widehat{u}_{i,h},\widehat{u}_{j,h})\big]_{\scriptscriptstyle 1\leq i\leq k, 1\leq j\leq k} \in\mathcal{R}^{k\times k}, \\
&\mathbf{b}_h =\big[b(\phi_{i,H},\widehat{u}_{j,h})\big]_{\scriptscriptstyle 1\leq i\leq N_H, 1\leq j\leq k} \in\mathcal{R}^{N_H\times k}, \\
&\boldsymbol{\beta} =\big[b(\widehat{u}_{i,h},\widehat{u}_{j,h})\big]_{\scriptscriptstyle 1\leq i\leq k, 1\leq j\leq k} \in\mathcal{R}^{k\times k}.
\end{split}
\end{equation*}
Equivalently, $\mathbf{a}_h$, $\boldsymbol{\alpha}$, $\mathbf{b}_h$ and $\boldsymbol{\beta}$ 
are also respectively given by $I_h^HA_h\widehat{\mathbf{u}}_h^{(\ell+1)}$, 
$(\widehat{\mathbf{u}}_h^{(\ell+1)})^TA_h\widehat{\mathbf{u}}_h^{(\ell+1)}$, $I_h^HB_h\widehat{\mathbf{u}}_h^{(\ell+1)}$ and $(\widehat{\mathbf{u}}_h^{(l+1)})^TB_h\widehat{\mathbf{u}}_h^{(\ell+1)}$.

Naturally, for the iterative initial guess of the eigenvalue problem \eqref{Augmented_GCGE}, we select the vector representation $\mathbf{u}_H=\mathbf{0}$ and $\gamma=\mathbf{I}$ of the current approximation function $\{\widehat{u}_{j,h}\}_{1\leq j\leq k}$ in the space $V_{H,h}$.

We would like to point out that the multi-vectors $\mathbf{a}_h$, $\mathbf{b}_h$ and $\mathbf{u}_H$ are stored in parallel, while $\boldsymbol{\alpha}$, $\boldsymbol{\beta}$ and $\gamma$ are stored sequentially with each process having a backup. We will later analyze this special matrix structure with both parallel and sequential storage, and explore how to accelerate computing this type of matrix.

\section{Precondition}
\label{sec_Precond} 
In this section, we discuss how to optimize PASE algorithm in a parallel computing environment. We analyze the scalability bottleneck of Algorithm \ref{Algorithm_PASE_Algebric} and identify that the evident issue lies in step 5, 
where solving the eigenvalue problem \eqref{Augmented_GCGE} involves 
numerous inner products of the vectors in the augmented subspace $V_{H,h}^{(\ell+1)}$. The eigenvalue problem  (\ref{Augmented_GCGE})
can be written as follows: Find $(\Lambda_h, [\mathbf{u}_H,\boldsymbol{\gamma}])$ such that
\begin{eqnarray}\label{Eigenvalue_Problem_AHh}
\left(
\begin{array}{cc}
A_H & \mathbf{a}_h\\
\mathbf{a}_h^\top  & \boldsymbol{\alpha}
\end{array}
\right)
\left(
\begin{array}{c}
\mathbf{u}_H\\
\boldsymbol{\gamma}
\end{array}
\right) =
\left(
\begin{array}{cc}
B_H & \mathbf{b}_h\\
\mathbf{b}_h^\top  & \boldsymbol{\beta}
\end{array}
\right)
\left(
\begin{array}{c}
\mathbf{u}_H\\
\boldsymbol{\gamma}
\end{array}
\right) \Lambda_h.
\end{eqnarray}

To clarify our motivation for designing a more scalable numerical method for the eigenvalue problem \eqref{Eigenvalue_Problem_AHh}, consider that the compound matrices $A_{Hh}$ and $B_{Hh}$ have a storage structure that differs from standard parallel matrix storage.
Based on the storage properties mentioned before of the matrices  $A_H$, $\mathbf{a}_h$, $\boldsymbol{\alpha}$, $B_H$, $\mathbf{b}_h$ and $\boldsymbol{\beta}$,
as an example, let us consider the computing details for the matrix-vector product $A_{Hh}\mathbf{u}_{Hh}$.
The matrix-vector product $A_{Hh}\mathbf{u}_{Hh}$ has the following form
\begin{eqnarray}\label{Matrix_Vector_Product}
A_{Hh}\mathbf{u}_{Hh} =
\left(
\begin{array}{cc}
A_H & \mathbf{a}_h\\
\mathbf{a}_h^\top  & \boldsymbol{\alpha}
\end{array}
\right)
\left(
\begin{array}{c}
\mathbf{u}_H\\
\boldsymbol{\gamma}
\end{array}
\right)
=
\left(
\begin{array}{c}
A_H\mathbf{u}_H+\mathbf{a}_h\boldsymbol{\gamma}\\
\mathbf{a}_h^\top \mathbf{u}_H+\boldsymbol{\alpha}\boldsymbol{\gamma}
\end{array}
\right).
\end{eqnarray}
The computation for matrix-vector product (\ref{Matrix_Vector_Product}) 
in the parallel computing environment has the following properties:
\begin{itemize}
\item The sparse matrix-vector product $A_H\mathbf{u}_H$ is efficient
since it only requires local communications.

\item The linear combination of vectors $\mathbf{a}_h\boldsymbol{\gamma}$ 
has no communication and can be accelerated by using BLAS3.

\item The inner product $\mathbf{a}_h^\top \mathbf{u}_H$ of two multi-vectors can also be accelerated by using BLAS3.
 Unfortunately, this computation involves global communication, 
 which will be the primary constraint on scalability.

\item The matrix-matrix product $\boldsymbol{\alpha}\boldsymbol{\gamma}$ of two dense 
matrices $\boldsymbol{\alpha}$ and $\boldsymbol{\gamma}$ has no communication
and can be efficiently implemented by using BLAS3.
\end{itemize}

Since the computation $\mathbf{a}_h^\top \mathbf{u}_H$ presents a scalability challenge, 
we design a type of preprocessing technique to avoid this part of computation. 
Since we have the decomposition for the matrix $A_{Hh}$
\begin{eqnarray}\label{Transformation_0}
\begin{pmatrix}
A_H& \mathbf{a}_h\\
\mathbf{a}_h^\top  &  \boldsymbol{\alpha}
\end{pmatrix}
=
\begin{pmatrix}
I_H & O \\
\mathbf{a}_h^{\top}A_H^{-1} & I_m
\end{pmatrix}
\begin{pmatrix}
A_H & O \\
O & \boldsymbol{\alpha} - \mathbf{a}_h^{\top} A_H^{-1} \mathbf{a}_h
\end{pmatrix}
\begin{pmatrix}
I_H & A_H^{-1} \mathbf{a}_h \\
O & I_m
\end{pmatrix},
\end{eqnarray}
the eigenvalue problem \eqref{Eigenvalue_Problem_AHh}
can be transformed to the following equivalent form
\begin{eqnarray}\label{Eigen_Problem_Hh_2}
\begin{pmatrix}
A_H & O \\
O & \widetilde{ \boldsymbol{\alpha}}
\end{pmatrix}
\begin{pmatrix}
\widetilde{\mathbf{u}}_H \\ \widetilde{\boldsymbol{\gamma}}
\end{pmatrix}
=
\begin{pmatrix}
B_H & \widetilde{\mathbf{b}_h}\\
\widetilde{\mathbf{b}}_h^{\top} & \widetilde{\boldsymbol{\beta}}
\end{pmatrix}
\begin{pmatrix}
\widetilde{\mathbf{u}}_H \\
\widetilde{\boldsymbol{\gamma}}
\end{pmatrix}\Lambda_h,
\end{eqnarray}
where
\begin{eqnarray*}
&&\left(
\begin{array}{c}
\widetilde{\mathbf{u}}_H \\
\widetilde{\boldsymbol{\gamma}}
\end{array}
\right)
=
\left(
\begin{array}{cc}
I_H & A_H^{-1} \mathbf{a}_h \\
O & I_m
\end{array}
\right)
\left(
\begin{array}{c}
\mathbf{u}_H \\\boldsymbol{\gamma}
\end{array}
\right),\\
&&
\widetilde{\boldsymbol{\alpha}} =  \boldsymbol{\alpha} - \mathbf{a}_h^{\top} A_H^{-1}\mathbf{a}_h,\\
&&
\widetilde{\mathbf{b}}_h = \mathbf{b}_h - B_HA_H^{-1}\mathbf{a}_h,\\
&&
\widetilde{\boldsymbol{\beta}} = \boldsymbol{\beta} - \mathbf{b}_h^{\top}A_H^{-1}\mathbf{a}_h - \mathbf{a}_h^{\top} A_H^{-1}\mathbf{b}_h
+ (A_H^{-1} \mathbf{a}_h)^{\top}B_H A_H^{-1} \mathbf{a}_h.
\end{eqnarray*}
In the following text, we will refer to this preprocessing method as \text{PRECOND-A}.
The advantage of solving the problem \eqref{Eigen_Problem_Hh_2} is that 
it avoids the scalability issue associated with the computation $\mathbf{a}_h^\top u_H$.

Given that the algorithm involves a significant number of orthogonalization operations concerning $B_{Hh}$, 
we also consider applying a preconditioner to $B_{Hh}$ beforehand. 
Similarly, $B_{Hh}$ has the following decomposition:
\begin{eqnarray*}
\begin{pmatrix}
B_H& \mathbf{b}_h\\
\mathbf{b}_h^\top  &  \boldsymbol{\beta}
\end{pmatrix}
=
\begin{pmatrix}
I_H & O \\
\mathbf{b}_h^{\top}B_H^{-1} & I_m
\end{pmatrix}
\begin{pmatrix}
B_H & O \\
O & \boldsymbol{\beta} - \mathbf{b}_h^{\top} B_H^{-1} \mathbf{b}_h
\end{pmatrix}
\begin{pmatrix}
I_H & B_H^{-1} \mathbf{b}_h \\
O & I_m
\end{pmatrix},
\end{eqnarray*}
The next steps are similar to those mentioned above and we will refer to this preprocessing method as \text{PRECOND-B}.
Based on the above transformation applied to $B_{Hh}$, when solving the augmented subspace problem, we can further improve the parallel efficiency of the linear system solving step by applying 
a transformation similar to \eqref{Transformation_0}  to the coefficient matrix during the linear system solution process. 
In the upcoming sections, we will denote this preprocessing method as \text{PRECOND-B-A}.


\begin{algorithm}[htb!]
\caption{GCG\_AUG algorithm}
\label{GCGAUG_Algorithm}
\begin{algorithmic}[1]
\State Choose ${\tt nev}$ vectors to build the block $X$ and two null blocks 
$P=[\ ]$, $W=[\ ]$.
Define the following small scale eigenvalue problem:
\begin{eqnarray*}
X^\top A_{Hh}XC=X^\top  B_{Hh}XC\Lambda.
\end{eqnarray*}
Solve this eigenvalue problem and choose the desired eigenvalues $\Lambda$ 
and eigenvectors $C$.
Then update $X\leftarrow XC$.

\State Generating $W$ by solving linear equations $W = (A_{Hh}-\mu B_{Hh})^{-1}(B_{Hh}X\Lambda)$ 
with the linear solver
and initial value $X$, where the shift $\mu$ is selected dynamically such that 
$(A_{Hh}-\mu B_{Hh})$ is a
symmetric positive definite matrix.

\State Define $V=[X,P,W]$ and do orthogonalization to $V$ in the sense of 
inner product deduced by the matrix $B_{Hh}$.

\State Solve the Rayleigh-Ritz problem $V^\top  A_{Hh}VC= C\Lambda$ to
obtain new desired eigenpair approximations $(\Lambda, C)$.
Then update $X_{\tt  new}= VC$.

\State Check the convergence of $\Lambda$ and $X_{\tt  new}$.  
If they converged, the iteration stops.

\State Otherwise, compute $P = X_{\tt  new} \backslash X$ and update $X = X_{\tt  new}$.
Then go to Step 2 for the next iteration until convergence.
\end{algorithmic}
\end{algorithm}

Considering the complexity of matrix and vector structures,
we need to employ an eigensolver that is both matrix-free and vector-free. Based on this consideration, we have made improvements to the GCG algorithm. 
Assume we need to compute the smallest {\tt nev} eigenpairs of \eqref{Eigen_Problem_Hh_2},
the GCG\_AUG algorithm can be defined by Algorithm \ref{GCGAUG_Algorithm}.
A series of optimization techniques have already been proposed in GCGE, such as an efficient block orthogonalization method, 
parallel computing for Rayleigh-Ritz problems, moving mechanics, 
and block conjugate gradient iteration to improve the efficiency and 
scalability for computing plenty of eigenpairs of large scale eigenvalue problems.
These techniques can also be applied to GCG\_AUG algorithm.
For more information about GCGE, please refer to 
\cite{LiWangXie,LiXieXuYouZhang,ZhangLiXieXuYou}.

\begin{algorithm}[htb!]
\caption{Solving eigenvalue problem \eqref{Eigenvalue_Problem_AHh} for parallel augmented subspace method}
\label{Algorithm_AHh}
\begin{algorithmic}[1]
\State Compute $\widehat{\mathbf{a}_h} = A_H^{-1} \mathbf{a}_h$ and then compute $\widetilde{\boldsymbol{\alpha}}$, $\widetilde{\boldsymbol{\beta}}$ and $\widetilde{\mathbf{b}_h}$ for (\ref{Eigen_Problem_Hh_2}) as follows:
\[
\widehat{\boldsymbol{\alpha}} = \boldsymbol{\alpha} - \widehat{ \mathbf{a}_h}^{\top}\mathbf{a}_h, \quad
\widetilde{\mathbf{b}_h} = \mathbf{b}_h - B_H\widehat{\mathbf{a}_h},
\]
\[
\widetilde{\boldsymbol{\beta}} = \boldsymbol{\beta} - \mathbf{b}_h^\top \widehat{\mathbf{a}_h} - \widehat{\mathbf{a}_h}^{\top} \mathbf{b}_h + \widehat{\mathbf{a}_h}^{\top}B_H\widehat{\mathbf{a}_h}.
\]

\State Solve the eigenvalue problem \eqref{Eigen_Problem_Hh_2} with GCG\_AUG defined by Algorithm \ref{GCGAUG_Algorithm} to get the eigenpair $(\Lambda_h, [\widetilde{u_H}, \widetilde{\boldsymbol{\gamma}}])$.

\State Perform the following transformation to get the eigenpair approximation
$(\Lambda, [u_H,\boldsymbol{\gamma}])$
\begin{eqnarray*}
&&
\left(
\begin{array}{c}
u_H \\\boldsymbol{\gamma}
\end{array}
\right)=
\left(
\begin{array}{cc}
I_H & -\widehat{\mathbf{a}_h} \\
O & I_m
\end{array}
\right)
\left(
\begin{array}{c}
\widetilde{u}_H \\
\widetilde{\gamma}
\end{array}
\right).
\end{eqnarray*}
\end{algorithmic}
\end{algorithm}

Based on the above discussion, 
the parallel numerical scheme for solving the eigenvalue problem \eqref{Aug_Eigenvalue_Problem_k} with method \text{PRECOND-A}
is outlined in Algorithm \ref{Algorithm_AHh}. The methods \text{PRECOND-B} and \text{PRECOND-B-A} are similar.
In this algorithm, the transformation from (\ref{Eigenvalue_Problem_AHh}) to  (\ref{Eigen_Problem_Hh_2}) can reduce much communication and improve
the scalability when solving the linear equations in step 2 of Algorithm \ref{GCGAUG_Algorithm}.


\section{Batch scheme for computing plenty of eigenpairs}
Since $\boldsymbol{\alpha}$, $\boldsymbol{\beta}$, and $\boldsymbol{\gamma}$ in \eqref{Eigenvalue_Problem_AHh} are stored in memory as serial dense matrices, with each process maintaining identical copies, and the size of these matrices, as in Algorithm \ref{Algorithm_PASE_Algebric}, is the number of desired eigenpairs.
The significant memory consumption and the reduced scalability when computing a large number of eigenpairs can impact the efficiency of the algorithm. Additionally, as mentioned in Section \ref{sec_Precond}, the dense vector blocks $\mathbf{a}_h$ and $\mathbf{b}_h$ also require substantial memory when the number of desired eigenpairs is large.
In order to improve efficiency further and reduce memory consumption 
for computing plenty of eigenpairs, we design a type of batch scheme inspired by the parallel methods mentioned in \cite{XuXieZhang}. With the help of Algorithms \ref{Algorithm_PASE_Algebric}, we come to present the batch scheme to compute the first $m(m \gg 1)$ eigenpairs of the eigenvalue problem \eqref{Weak_Eigenvalue_Discrete}. The batch scheme here is a type of eigenpair-wise parallel method for eigenvalue problems.






Specifically, we partition the $m$ desired eigenpairs into $\tau$ batches and solve them sequentially. 
Let the $i$-th batch of eigenpairs be denoted by $\lambda_{i,h}$, $i=m_i,...,m_i+k_i-1$, 
and let the current approximation of eigenfunctions in $V_h$ be $\widetilde{u}_{i,h}$ for $i=m_i,...,m_i+k_i-1$. 
We then solve the following problem in the augmented subspace $V_{H,h}=V_H+{\rm span}\{\widetilde{u}_{m_i,h}, ..., \widetilde{u}_{m_i+k_i-1,h}\}$:
\begin{equation}\label{Batch_Aug_FEM}
  a(u_{i,h},v_{H,h}) = \lambda_{i,h}b(u_{i,h},v_{H,h}),
  \ \  \forall v_{H,h}\in V_{H,h},\ \  i=m_i, \cdots, m_i+k_i-1.
\end{equation}

The challenge with this computational approach is that, starting from the second batch, 
the required eigenpairs are no longer the smallest eigenpairs of problem \eqref{Batch_Aug_FEM}, 
but instead a subset of internal eigenpairs. 
For the $i$-th batch, $i=1,...,\tau$, corresponding to eigenpairs from the $m_i$-th to the $(m_i+k_i-1)$-th, we first estimate value $\theta_i=(\lambda_{m_i, H}+\lambda_{m_i+k_i-1, H})/2$,
where the eigenvalues $\lambda_{1,H},...,\lambda_{m, H}$ are obtained by solving \eqref{initial_eigenproblem} on the coarse space.
We then solve problem \eqref{Batch_Aug_FEM} to find the $n$ eigenvalues closest to $\theta$,
along with their corresponding eigenvectors.

It is important to note that the $k_i$ eigenpairs needed for solving the $i$-th batch are not necessarily the $k_i$ eigenpairs closest to $\theta_i$ of \eqref{Batch_Aug_FEM}. 
In parctice, we often need to ensure that the computed set of eigenpairs includes more than $k_i$ pairs, say $n>k_i$, to guarantee that the required ones are captured. 
Once the $n$ eigenpairs closest to $\theta$ are obtained, we select the $k_i$ eigenpairs that are most aligned with the current eigenspace. 
Specifically, we identify the $k_i$ eigenpairs from the $n$ computed ones 
whose eigenfunctions have the largest components in the space ${\rm span}\{\widetilde{u}_{m_i,h}, ..., \widetilde{u}_{m_i+k_i-1,h}\}$. 
These selected eigenpairs are then used as the solutions for the current iteration.

Assume the $n$ computed eigenpairs are $\{\lambda_{i,h},v_{i,h}\}$, $i=1,...n$, where
\begin{equation*}
    v_{i,h}=u_{i,H}+\sum_{j=0}^{k_i-1}\gamma_{j+1,i}\widetilde{u}_{m_i+j,h},\quad i=1,...,n
\end{equation*}
Define the projection of $v_{i,h}$ onto ${\rm span}\{\widetilde{u}_{m_i,h},...,\widetilde{u}_{m_i+k_i-1,h}\}$ in the sense of $b$-product as $$\Pi_hv_{i,h}=\sum_{j=0}^{k_i-1}x_{i,j+1}\widetilde{u}_{m_i+j,h}, \\\\ i=1,...,n.$$ 

Define $\mathbf X=(x_{i,j})$ for $i=1,\cdots,n$ and $j=1,\cdots,k_i$, we can derive the equality $\mathbf X\boldsymbol\beta = \mathbf{u}^T_H\mathbf{b}_h+\boldsymbol\gamma^T\boldsymbol\beta$ where $\boldsymbol\beta$ and $\mathbf{b}_H$ is defined by \eqref{Definition_Hh} and $\boldsymbol\gamma$ and $\mathbf{b}_H$ is obtained by \eqref{Augmented_GCGE}. 
Solving for $\mathbf X$ involves dense matrix operations. 
The components of $[v_{1,h}, \cdots, v_{n,h}]$ in the space $\text{span}\{\widetilde{u}_{m_i,h},\cdots,\widetilde{u}_{m_i+k_i-1,h}\}$ in the sense of $b$-product are given by
\begin{eqnarray}\label{largest_component}
\|\Pi_h v_{i,h}\|_b^2=X(i,:)\boldsymbol{\beta} X(i,:)^\top, \ \ \ \ 	i=1,...,n. 
\end{eqnarray}
Based on these component values, we select the $k_i$ eigenpairs with the largest components as the desired ones.

The next step is to efficiently and stably compute the $n$ eigenpairs near a specific shift $\theta$. We should modified the eigenvalue (\ref{Eigenvalue_Problem_AHh}) to the following form:
Find $(\Lambda_h, [u_H,\boldsymbol{\gamma}])$ such that
\begin{eqnarray}\label{Eigenvalue_Problem_AHh_shift}
\left(
\begin{array}{cc}
A_H-\theta B_H & \mathbf{a}_h-\theta \mathbf{b}_h\\
(\mathbf{a}_h-\theta \mathbf{b}_h)^\top  & \boldsymbol{\alpha}-\theta \boldsymbol{\beta}
\end{array}
\right)
\left(
\begin{array}{c}
\mathbf{u}_H\\
\boldsymbol{\gamma}
\end{array}
\right) =
\left(
\begin{array}{cc}
B_H & \mathbf{b}_h\\
\mathbf{b}_h^\top  & \boldsymbol{\beta}
\end{array}
\right)
\left(
\begin{array}{c}
\mathbf{u}_H\\
\boldsymbol{\gamma}
\end{array}
\right) \Lambda_h.
\end{eqnarray}
We can use the idea in Section \ref{sec_Precond} to improve the scalability of the batch scheme 
for computing plenty of eigenpairs. 
Based on (\ref{Transformation_0}), the eigenvalue problem \eqref{Eigenvalue_Problem_AHh_shift}
can also be transformed to the following equivalent form
\begin{eqnarray}\label{Modified_Eigen_Problem_Hh_2}
\begin{pmatrix}
\widehat A_H & O \\
O & \widetilde{\boldsymbol{\alpha}}
\end{pmatrix}
\begin{pmatrix}
\widetilde{\mathbf{u}}_H \\ \widetilde{\boldsymbol{\gamma}}
\end{pmatrix}
=
\begin{pmatrix}
B_H & \widetilde{\mathbf{b}}_h \\
\widetilde{\mathbf{b}}_h^{\top} & \widetilde{\boldsymbol{\beta}}
\end{pmatrix}
\begin{pmatrix}
\widetilde{\mathbf{u}}_H \\
\widetilde{\boldsymbol{\gamma}}
\end{pmatrix}\Lambda_h,
\end{eqnarray}
In order to compute the correct desired eigenpairs,
we need to modify GCG\_AUG algorithm with a shift $\theta$ such that it can compute 
the interior eigenpairs.  Assume we need to compute the closest 
{\tt nev} eigenpairs to $\theta$, the modified GCG\_AUG algorithm can be
defined by Algorithm \ref{Modified_GCG_Algorithm}.
\begin{algorithm}[htb!]
\caption{Modified GCG\_AUG algorithm for (\ref{Modified_Eigen_Problem_Hh_2}) with shift $\theta$}
\label{Modified_GCG_Algorithm}
\begin{algorithmic}[1]
\State Update $A_{Hh}$ and $B_{Hh}$ with shift $\theta$: $A_{Hh} \leftarrow A_{Hh} - \theta B_{Hh}$.

\State Choose ${\tt nev}$ vectors to build the block $X$ and two null blocks $P=[\ ]$, $W=[\ ]$. 
Define the following small scale eigenvalue problem:
\begin{eqnarray*}
X^\top A_{Hh}^\top A_{Hh}XC=X^\top A_{Hh}^\top B_{Hh}XC\Lambda.
\end{eqnarray*}
Solve this eigenvalue problem and choose the desired eigenvalues $\Lambda$ and eigenvectors $C$. 
Then update $X \leftarrow XC$.

\State Generate $W$ by solving the linear equations $W = A_{Hh}^{-1}(B_{Hh} X \Lambda)$ with the linear solver and the initial value $X$.

\State Define $V = [X, P, W]$ and perform orthogonalization to $V$ in the sense of $L^2$ inner product.

\State Solve the harmonic Rayleigh-Ritz problem:
\[
V^\top A_{Hh}^\top A_{Hh} V C = V^\top A_{Hh}^\top B_{Hh} V C \Lambda,
\]
to obtain new desired eigenpair approximations $(\Lambda, C)$. Then update $X_{\tt new}$ as $X_{\tt new} = VC$.

\State Check the convergence of $\Lambda$ and $X_{\tt new}$. If they have converged, the iteration stops.

\State Otherwise, compute $P = X_{\tt new} \backslash X$ and update $X = X_{\tt new}$. Then go to Step 2 for the next iteration until convergence.

\State Solve the Rayleigh-Ritz problem:
\[
V^\top A_{Hh} V C = V^\top B_{Hh} V C \Lambda,
\]
to obtain the desired eigenpair approximations $(\Lambda, C)$.
\end{algorithmic}
\end{algorithm}


\begin{algorithm}[htbp!]
\caption{Algebraic version of PASE with Batch scheme}
\label{Algorithm_PASE_batch}
\begin{algorithmic}[1]
\State Solve the initial small eigenvalue problem for the first $m$ eigenpairs:
\[
A_H \mathbf{u}_H = B_H \mathbf{u}_H \Lambda_H.
\]

\State Batch the eigenpairs into $\tau$ batches, where the $i$-th batch includes the $m_i$-th to ($m_i + k_i - 1$)-th eigenpairs. Set $i = 0$.

\State Do the following iterations:
\begin{enumerate}[(a)]
    \item Set $\ell = 1$ and interpolate the multi-vectors $\mathbf{u}_{i,H} = \mathbf{u}_H[m_i:m_i+k_i-1]$ to obtain $\mathbf{u}_{i,h}^{(\ell)} = I_H^h \mathbf{u}_{i,H}$ and define $\Lambda_{i,h}^{(\ell)} = \Lambda_H[m_i:m_i+k_i-1]$.
    
    \item Define the shift $\theta = (\lambda_{m_i,H} - \lambda_{m_i+k_i-1,H}) / 2$.

    \item Solve the following linear equation by performing $s$ CG iteration steps with the initial solution $\mathbf{u}_{i,h}^{(\ell)}$ to obtain $\widehat{\mathbf{u}}_{i,h}^{(\ell+1)}$:
    \[
    A_h \widehat{\mathbf{u}}_{i,h}^{(\ell+1)} = B_h \mathbf{u}_{i,h}^{(\ell)} \Lambda_{i,h}^{(\ell)}.
    \]

    \item Assemble the composite matrices $A_{Hh}$ and $B_{Hh}$ as \eqref{Definition_Hh} and solve the following eigenvalue problem for the $n > k_i$ eigenpairs closest to the shift $\theta$ such that they contain the desired ones:
    \[
    A_{Hh} \mathbf{u}_{Hh} = B_{Hh} \mathbf{u}_{Hh} \Lambda_h, \quad \mathbf{u}_{Hh} \triangleq \left( \begin{array}{c} \mathbf{u}_H \\ \boldsymbol{\gamma} \end{array} \right).
    \]
    
    \item Choose $k$ eigenpairs which are closest to the current eigenspace by \eqref{largest_component}. Define them as:
    \[
    \mathbf{u}_{i,Hh} \triangleq \left( \begin{array}{c} \mathbf{u}_{i,H} \\ \boldsymbol{\gamma}_i \end{array} \right)
    \]
    and compute:
    \[
    \widehat{\mathbf{u}}_{i,h}^{(\ell+1)} = I_H^h \mathbf{u}_{i,H}^{(\ell+1)} + \boldsymbol{\gamma}_i \widehat{\mathbf{u}}_{i,h}^{(\ell+1)}.
    \]

    \item Solve the following linear equation by performing $s$ CG iteration steps with the initial solution $\widehat{\mathbf{u}}_{i,h}^{(\ell+1)}$ to obtain ${\mathbf{u}}_{i,h}^{(\ell+1)}$:
    \[
    A_h \mathbf{u}_{i,h}^{(\ell+1)} = B_h \widehat{\mathbf{u}}_{i,h}^{(\ell+1)} \Lambda_{i,h}^{(\ell+1)}.
    \]

    \item Set $\ell = \ell + 1$ and go to Step (c) for the next iteration until convergence.
\end{enumerate}

\State Set $i = i + 1$ and go to Step 3 for the next batch until $i > \tau$.
\end{algorithmic}
\end{algorithm}

Using Algorithms \ref{Modified_GCG_Algorithm}, 
we propose a batch scheme for computing the first $m$ eigenpairs of the eigenvalue problem \eqref{Weak_Eigenvalue_Discrete} with the augmented subspace method. 
The detailed algebraic scheme is outlined in Algorithm \ref{Algorithm_PASE_batch}. When computing plenty of eigenpairs, 
whether and how to perform the batch scheme mainly depend on the number of required eigenpairs.


\section{Numerical results}\label{Section_Numerical_Results}
In this section, we will sequentially demonstrate the computational efficiency of PASE through several examples. We will highlight how our techniques significantly enhance algorithm performance and their effectiveness in problems that require adaptive refinement methods.

In numerical experiments, we adopt the block conjugate gradient (BCG) iteration method with a maximum of $40$ steps to solve the linear problem \eqref{presoomthing}
 and \eqref{postsoomthing} in Algorithm \ref{Algorithm_PASE_Algebric}. In order to show the efficiency clearly, we also compare the numerical results of
the proposed method with that of Krylov-Schur and LOBPCG methods from SLEPc \cite{SLEPC}. Furthermore, all the methods use the same convergence criterion which is set to be $\|Ax-\lambda Bx\|_2/|\lambda|\leq 1{\tt e}$-$8$ for the algebraic eigenvalue problem: $Ax=\lambda Bx$, where $\|\cdot\|_2$ denotes the $L^2$-norm for vectors.
 
The numerical examples are carried out on LSSC-IV in the State Key Laboratory of 
Scientific and Engineering
Computing, Chinese Academy of Sciences. Each computing node has two $18$-core 
Intel Xeon Gold $6140$
processors at $2.3$ GHz and $192$ GB memory. The matrices are produced by 
the open package of open parallel finite element method(OpenPFEM). 
OpenPFEM provides the parallel finite element discretization for the
partial differential equations and corresponding eigenvalue problems. 
The source code can be downloaded from \url{https://gitlab.com/xiegroup/pase1.0}.

\subsection{Tests for the model eigenvalue problem}

In this subsection, we examine a typical eigenvalue problem that demonstrates the high computational efficiency of PASE: Find $(\lambda,u)\in\mathcal R\times V$ such that $||u||_a=1$ and

\begin{eqnarray}\label{AMG2D}
\left\{
\begin{array}{rcl}
-\Delta u &=&\lambda u,\ \ \ {\rm in}\ \Omega,\\
u&=&0,\ \ \ \ \ {\rm on}\ \partial\Omega,
\end{array}
\right.
\end{eqnarray}
where $\Omega=(0,1)^2$.

From the perspective of algorithm implementation, the process of using AMG (algebraic multigrid) method and GMG (geometric multigrid) method in Algorithm \ref{Algorithm_PASE_Algebric} is the same. In fact, the abstract multigrid constructed through the AMG method can be applied to solve algebraic eigenvalue problems in cases where geometric multigrid information is in sufficient, even though there is no true geometric grid. 
Therefore, we will solve this problem using both methods, where the AMG method is based on the parallel implementation BoomerAMG from hypre \cite{hypre} and the GMG method is generated by OpenPFEM.

We use Algorithm \ref{Algorithm_PASE_Algebric} with both AMG and GMG methods to compute the first $200$ and $800$ eigenpairs for problem \eqref{AMG2D}, comparing the results with GCGE, LOBPCG and Krylov-Schur methods. The LOBPCG and Krylov-Schur solvers used are the built-in versions from SLEPc. When performing comparison tests, we have already set optimal parameters for both solvers, and this approach is applied to all subsequent test cases. 
When comparing with PASE, these solvers solve the eigenvalue problem corresponding to the finest level mesh.
The CPU times of computation are listed in Tables \ref{table1} and \ref{table2}. We use "nDofs" to represent the number of degrees of freedom, and this abbreviation is also used in the following tables.
The results listed in the tables demonstrate the advantages of PASE in solving eigenvalue problems, particularly for large-scale eigenvalue problems. 


\begin{table}[!hbt]
\caption{The CPU time (in seconds) of Algorithm \ref{Algorithm_PASE_Algebric}, GCGE, LOBPCG and Krylov-Schur 
for computing the first $200$ eigenpairs of (\ref{AMG2D}).}\label{table1}
\begin{tabular}{|c|c|c|c|c|c|c|}\hline
nDofs     & nproc & PASE(AMG) & PASE(GMG) & GCGE & LOBPCG & Krylov-Schur \\ \hline
261,121   &  12   &   30.21   &  26.83    & 39.09     & 242.35 &  57.18     \\ \hline
1,046,529 &  36   &   86.71   &  47.33    & 91.05     & 936.82 & 116.68     \\ \hline
4,190,209 &  108  &   236.35  &  53.81    & 301.59    & 1125.90& 315.39  \\ \hline
\end{tabular}
\end{table}

\begin{table}[!hbt]
\caption{The CPU time (in seconds) of Algorithm \ref{Algorithm_PASE_Algebric}, GCGE, LOBPCG and Krylov-Schur 
for computing the first $800$ eigenpairs of (\ref{AMG2D}),
where the symbol $``-"$ means the computer runs out of memory.}\label{table2}
\begin{tabular}{|c|c|c|c|c|c|c|}\hline
nDofs   & nproc & PASE(AMG) & PASE(GMG) & GCGE & LOBPCG & Krylov-Schur \\ \hline
1,046,529 &  36   & 425.02    & 528.33    & 369.32   &5828.60 &496.92    \\ \hline
4,190,209 &  108  & 674.26    & 557.89    & 1034.43  &11744.66 &1295.94     \\ \hline
16,769,025&  324  & 910.53    & 894.72    & 3564.73   &12053.51 &3282.26    \\ \hline
\end{tabular}
\end{table}

However, it is important to mention that, due to the lack of precise projection operators and restriction operators during AMG method, it often requires more iterations than geometric multigrid to achieve the same level of accuracy. Moreover, as discussed in the theoretical analysis in Section \ref{sec_2}, the convergence speed of PASE is also impacted by the effectiveness of the coarse grid in approximating the low-frequency information. 

In addition, our approach involves applying AMG method to the stiffness matrix $A_h$ to obatin the coarse one $A_H$ and the interpolation operator $I_H^h$ and restriction operator $I_h^H$, followed by the computation $B_H=I_h^HB_hI_H^h$. This approach may also compromise the quality of the coarse space.

\subsection{Tests for Precondition}
\label{sec_mep}
In this subsection, we use Algorithm \ref{Algorithm_PASE_Algebric} with different precondition methods to solve the 
following model eigenvalue problem:
Find $(\lambda,u)\in\mathcal R\times V$ such that $\|u\|_a=1$ and
\begin{eqnarray}\label{Model_Eigenvalue_Problem}
\left\{
\begin{array}{rcl}
-\Delta u &=&\lambda u,\ \ \ {\rm in}\ \Omega,\\
u&=&0,\ \ \ \ \ {\rm on}\ \partial\Omega,
\end{array}
\right.
\end{eqnarray}
where $\Omega = (0,1)^3$.

\begin{figure}[!htb]
    \centering
    \includegraphics[height=5.5cm,width=6cm]{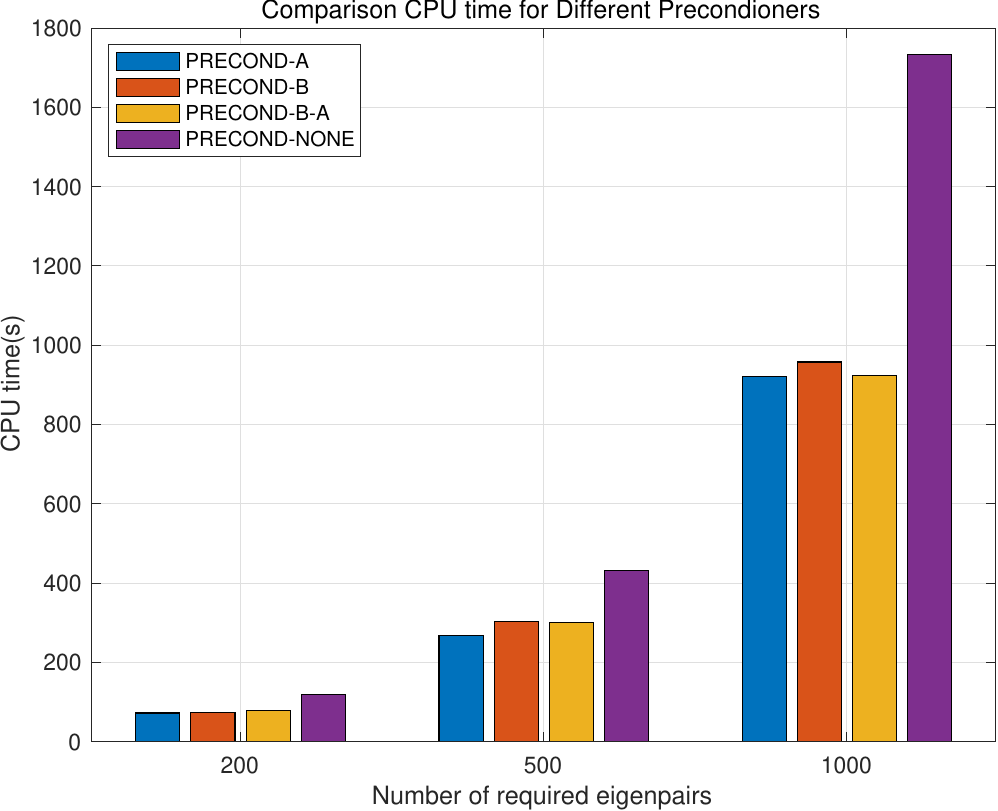}
    \caption{The CPU time of Algorithm \ref{Algorithm_PASE_Algebric} with different preconditioners 
    for example 2 with $\rm nproc = 288$ and $\rm tol = 10^{-8}$.}
    \label{fig:precond_ex1}
    \end{figure}

We first check the computational time of PASE for computing 
the first $200$, $500$, and $1000$ eigenpairs of problem \eqref{Model_Eigenvalue_Problem} 
on the finest level of finite element space with the number of degree of freedoms 
$13,785,215$ with several different preconditioning methods.
 The results presented in Figure \ref{fig:precond_ex1} indicate that \text{PRECOND-A} exhibits the highest computational efficiency. 
 Furthermore, as the number of eigenvalue pairs to be computed increases, the benefits of the preconditioner become more pronounced. When $\text{nev} = 1000$, using \text{PRECOND-A} can nearly double the computational efficiency compared to not employing a preconditioner. Therefore, in the tests conducted in this section, we default to using \text{PRECOND-A}.

    And we also investigate the
    efficiency and scalability for computing the first $200$ and $800$ eigenpairs, respectively.
    Tables \ref{table3} and \ref{table4} show the corresponding numerical results for computing 
    the first $200$ and $800$ eigenpairs by Algorithm \ref{Algorithm_PASE_Algebric}, GCGE, Krylov-Schur 
    and LOBPCG method. The results presented in Tables \ref{table3} and \ref{table4} 
    reveal that PASE exhibits superior computational efficiency to other eigensolvers. 
    Moreover, as the number of degrees of freedom increases, 
    the efficiency advantage of PASE becomes increasingly pronounced.

    \begin{table}[!htb]
    \caption{The CPU time (in seconds) of Algorithm \ref{Algorithm_PASE_Algebric}, GCGE, LOBPCG and Krylov-Schur 
    for computing the first $200$ eigenpairs of (\ref{Model_Eigenvalue_Problem}), where
    the symbol ``$-$" means the computer runs out of memory.}\label{table3}
    \begin{tabular}{|c|c|c|c|c|c|}\hline
     nDofs  & nproc & Time of PASE & Time of GCGE & Time of LOBPCG & Time of KS \\ \hline
    1,698,879       & 72    & 58.25        & 112.00       & 516.83         & 342.767 \\ \hline
    13,785,215      &288    & 71.06        & 415.88       & 2244.152       & $-$     \\ \hline
    111,063,295     &576    & 443.67       & 4359.75      & $-$            & $-$     \\ \hline
    \end{tabular}
    \end{table}
    \begin{table}[!htb]
    \caption{The CPU time (in seconds) of Algorithm \ref{Algorithm_PASE_Algebric}, GCGE, LOBPCG and Krylov-Schur 
    for computing the first $800$ eigenpairs of (\ref{Model_Eigenvalue_Problem}),
    where the symbol $``-"$ means the computer runs out of memory.}\label{table4}
    \begin{tabular}{|c|c|c|c|c|c|}\hline
    nDofs  & nproc & Time of PASE & Time of GCGE & Time of LOBPCG & Time of KS \\ \hline
    1,698,879       &72     & 625.69       & 845.07       & 2603.90        & $-$   \\ \hline
    13,785,215      &288    & 609.49       & 3101.74      & 9811.37        & $-$   \\ \hline
    111,063,295     &1440   & 843.58       & 7010.06      & $-$            & $-$   \\ \hline
    \end{tabular}
    \end{table}
    
\begin{figure}[!htb]
    \centering
    \includegraphics[height=5.5cm,width=5.5cm]{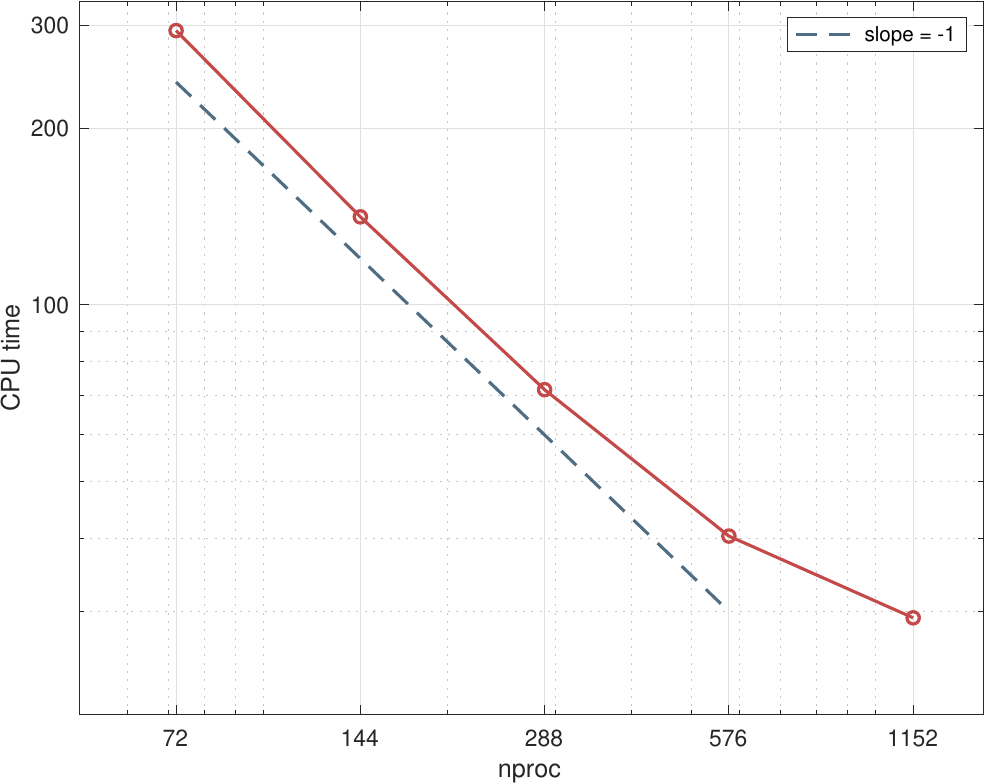}
    \caption{The Parallel Scalability of Algorithm \ref{Algorithm_PASE_Algebric} for example 2, 
    with $\rm nev = 200$ and $\rm tol = 10^{-8}$.}
    \label{fig:scaling_ex1}
    \end{figure}

    To assess the parallel scalability of algorithm \ref{Algorithm_PASE_Algebric}, 
    we perform experiments to evaluate the computation time of PASE for 
    calculating the first $200$ eigenvector pairs of the problem \eqref{Model_Eigenvalue_Problem} 
    on the finest level of finite element space with the number of degree of freedoms 
    $13,785,215$ under different numbers of processors. Figure \ref{fig:scaling_ex1} shows the 
    corresponding CPU time (in seconds) for different number of processors.  
    From Figure \ref{fig:scaling_ex1}, it is evident that PASE exhibits excellent scalability. For this example, since the dimension of $A_H$ is only $206,367$, when the number of processes exceeds $1,000$, each process owns very few degrees of freedom and the communication occupies a significant portion of the CPU time, resulting in a slight decrease in scalability.
    

    \subsection{Tests for batch scheme}\label{sec_mgep}
    We select the following second order elliptic eigenvalue problem to test the proposed Algorithm \ref{Algorithm_PASE_batch} with batch scheme. The aim of the numerical tests is to evaluate the performance of the method in computing a large number of eigenpairs for extremely large-scale matrices: 
    Find $(\lambda,u)\in\mathcal R\times V$ such that $\|u\|_a=1$ and
    \begin{eqnarray}\label{general}
    \left\{
    \begin{array}{rcl}
    -\nabla\cdot(\mathcal{A}\nabla u) +\varphi u&=&\lambda u\ \ \ {\rm in}\ \Omega,\\
    u&=&0\ \ \ \ \ {\rm on}\ \partial\Omega,\\
    \|u\|_{0,\Omega}&=&1,
    \end{array}
    \right.
    \end{eqnarray}
    where
    \begin{equation*}
    \mathcal{A}=\left(
    \begin{array}{ccc}
    1+(x_1-\frac{1}{2})^2 & (x_1-\frac{1}{2})(x_2-\frac{1}{2}) & (x_1-\frac{1}{2})(x_3-\frac{1}{2})\\
    (x_1-\frac{1}{2})(x_2-\frac{1}{2})&1+(x_2-\frac{1}{2})^2&(x_2-\frac{1}{2})(x_3-\frac{1}{2})\\
    (x_1-\frac{1}{2})(x_3-\frac{1}{2})&(x_2-\frac{1}{2})(x_3-\frac{1}{2})&1+(x_3-\frac{1}{2})^2
    \end{array}
    \right),
    \end{equation*}
    $\varphi=e^{(x_1-\frac{1}{2})(x_2-\frac{1}{2})(x_3-\frac{1}{2})}$ and
    $\Omega = (0,1)^3$.

    We first illustrate the improvement in memory achieved by the batch scheme. We compute the initial 2,000 eigenpairs the maximum memory usage for computing the first 2000 eigenvalues with Dofs = 13,785,215 and the first 3000 eigenvalues with Dofs = 111,063,295. For each case, we compared the peak memory usage between Algorithm \ref{Algorithm_PASE_Algebric} and Algorithm \ref{Algorithm_PASE_batch} with batch scheme in Table \ref{table5}.


\begin{table}[!hbt]
    \centering
    \caption{The peak memory usage of Algorithm \ref{Algorithm_PASE_Algebric} and \ref{Algorithm_PASE_batch} for 
    computing the first 2000 and 3000 eigenpairs of (\ref{general}), where
    the symbol ``$-$" means that the peak memory surpasses the total memory capacity.}\label{table5}
    \begin{tabular}{ccccc}
        \toprule
        \textbf{nDofs} & \textbf{nev} & \textbf{Batching} & \textbf{nprocs} & \textbf{Peak Memory of PASE} \\ \midrule
        13,785,215      & 2000 & No       & 288    & $-$                 \\ 
        13,785,215      & 2000 & No       & 360    & 1,749,783 MB        \\ 
        13,785,215      & 2000 & Yes      & 288    & 618,382 MB          \\ 
        13,785,215      & 2000 & Yes      & 360    & 642,644 MB          \\ \midrule
        111,063,295     & 3000 & No       & 2880   & $-$                 \\ 
        111,063,295     & 3000 & No       & 4500   & 19,966,051 MB       \\ 
        111,063,295     & 3000 & Yes      & 2880   & 8,428,079 MB        \\ 
        111,063,295     & 3000 & Yes      & 4500   & 10,673,858 MB       \\ 
        111,063,295     & 5000 & Yes      & 4500   & 13,505,930 MB       \\ 
        \bottomrule
    \end{tabular}
\end{table}


    Next in Table \ref{table6}, we present the the improvements in computation time achieved through batching. We also provide a comprehensive comparison of the results with GCGE in table \ref{table6}, while illustrating the significant advancements of PASE with batch scheme. The results clearly show that batching not only reduces the memory consumption but also enhances the overall computational efficiency of PASE. 
    

\begin{table}[!hbt]
    \centering
    \small 
    \setlength{\tabcolsep}{3pt} 
    \renewcommand{\arraystretch}{1.2} 
    \caption{The CPU time (in seconds) of Algorithms for computing the first $2000$, $3000$ and $6000$ eigenpairs of (\ref{general}), where the symbol ``$-$" means out of memory.}\label{table6}
    \begin{tabular}{p{2cm} p{1cm} p{1.3cm} p{2cm} p{2cm} p{2cm}}
        \toprule
        \textbf{nDofs} & 
        \textbf{nev} & 
        \textbf{nprocs} & 
        \textbf{Without Batching Time} & 
        \textbf{Batching Time} & 
        \textbf{GCGE Time} \\ 
        \midrule
        13,785,215 & 2000 & 288 & $-$ & 3006.61 & 5464.00 \\ 
        13,785,215 & 2000 & 360 & 3522.72 & 2571.93 & 4856.57 \\ 
        \midrule
        111,063,295 & 3000 & 2880 & $-$ & 4713.13 & 13727.36 \\ 
        111,063,295 & 3000 & 4500 & 6108.42 & 4235.51 & 8804.24 \\ 
        111,063,295 & 6000 & 4500 & $-$ & 9723.55 & 21467.15 \\ 
        \bottomrule
    \end{tabular}
\end{table}

    The results presented in Tables \ref{table3} and \ref{table4} indicate that, 
    with the implementation of the batch strategy, PASE efficiently 
    solves a significant number of eigenpairs. 
    Furthermore, as the matrix size expands, the computational efficiency advantage of PASE over GCGE becomes increasingly apparent.

\subsection{PASE with Adaptive Finite Element Method}

In this subsection, we aim to evaluate the performance of PASE when applied to problems with singularities. In order to keep the scalability of the 
    finite element discretization, we use the redistributing process 
    of package OpenPFEM to balance the distribution of the meshes. PASE is designed to be compatible with the adaptive mesh refinement and redistribution strategy, which ensures an optimal mesh quality for singular problems.

In all the following examples, adaptive mesh refinement is employed to handle the singularities inherent in the problems. The residual-based a posteriori error estimation is used to generate error indicators for each mesh element during the adaptive refinement process. The element residual $\mathcal{R}_T\left(\Lambda_h, \Phi_h\right)$ will be defined separately for each example below, while the jump residual $\mathcal{J}_e\left(\Phi_h\right)$ for the eigenpair approximation $\left(\Lambda_h, \Phi_h\right)$ remains the same and is defined as follows:

\begin{eqnarray}\label{posterior} \mathcal{J}_e\left(\Phi_h\right):=\left(\left.\frac{1}{2} \nabla u_{i, h}\right|_{T^{+}} \cdot v^{+}+\left.\frac{1}{2} \nabla u_{i, h}\right|_{T^{-}} \cdot v^{-}\right)_{i=1}^N, 
\quad \text { on } e \in \mathcal{E}_h, \end{eqnarray}

where $\mathcal{E}_h$ denotes the set of interior faces of $\mathcal{T}_h$, and $e$ is the common side of elements $T^{+}$ and $T^{-}$ with unit outward normals $v^{+}$ and $v^{-}$. The local error indicator $\eta_k^2\left(\Lambda_h, \Phi_h, T\right)$ is defined as follows:

\begin{eqnarray*}
\eta_k^2\left(\Lambda_h, \Phi_h, T\right):=h_T^2\left\|\mathcal{R}_T\left(\Lambda_h, 
\Phi_h\right)\right\|_{0, T}^2+\sum_{e \in \mathcal{E}_{h_k, e \subset \partial T}} h_e\left\|\mathcal{J}_e\left(\Phi_h\right)\right\|_{0, e}^2 .
\end{eqnarray*}

The error estimate $\eta_k^2\left(\Lambda_h, \Phi_h, \omega\right)$ for a subset $\omega \subset \Omega$ is defined as:

\begin{eqnarray}\label{eta_k_w}
\eta_k^2\left(\Lambda_h, \Phi_h, \omega\right)=\sum_{T \in \mathcal{T}_h, T \subset \omega} \eta_k^2\left(\Lambda_h, \Phi_h, T\right) .
\end{eqnarray}

We use the D\"{o}rfler's marking strategy to mark elements for local refinement based on the error indicator (\ref{eta_k_w}).

\subsubsection{Example: Hydrogen Atom}
In this example, we consider the following eigenvalue problem for the Hydrogen atom:
Find $(\lambda,u)\in\mathcal R\times V$ such that $\|u\|_a=1$ and
\begin{eqnarray}\label{hy}
-\frac{1}{2}\Delta u-\frac{1}{|x|}u=\lambda u, \quad {\rm in} \  \Omega,
\end{eqnarray}
where $\Omega=\mathcal R^3$. The eigenvalues of (\ref{hy}) are $\lambda_h=-\frac{1}{2n^2}$ with multiplicity $n^2$ for any positive integer $n$. As $n$ increases, the spectral gap narrows, and the multiplicity grows, which increases the difficulty of solving this eigenvalue problem. The goal of this example is to demonstrate that PASE can efficiently compute clustered eigenvalues and their eigenfunctions.

Since the eigenfunction decays exponentially, we set $\Omega = (-4,4)^3$ with a boundary condition $u=0$ on $\partial\Omega$. We use the adaptive refinement strategy to handle the singularities, coupling it with PASE.

The element residual for the Hydrogen atom example is defined as:
\begin{eqnarray}\label{residual_hydrogen}
\mathcal{R}_T\left(\Lambda_h, \Phi_h\right):=\left(\frac{1}{2}\Delta u_{i, h}
-\frac{1}{|x|}u_{i, h}+\lambda_{i, h} u_{i, h}\right)_{i=1}^N, 
\quad \text { in } T \in \mathcal{T}_{h}.
\end{eqnarray}

\begin{figure}[!htb]
\centering
\begin{minipage}{0.48\textwidth}
  \centering
  \includegraphics[width=\textwidth]{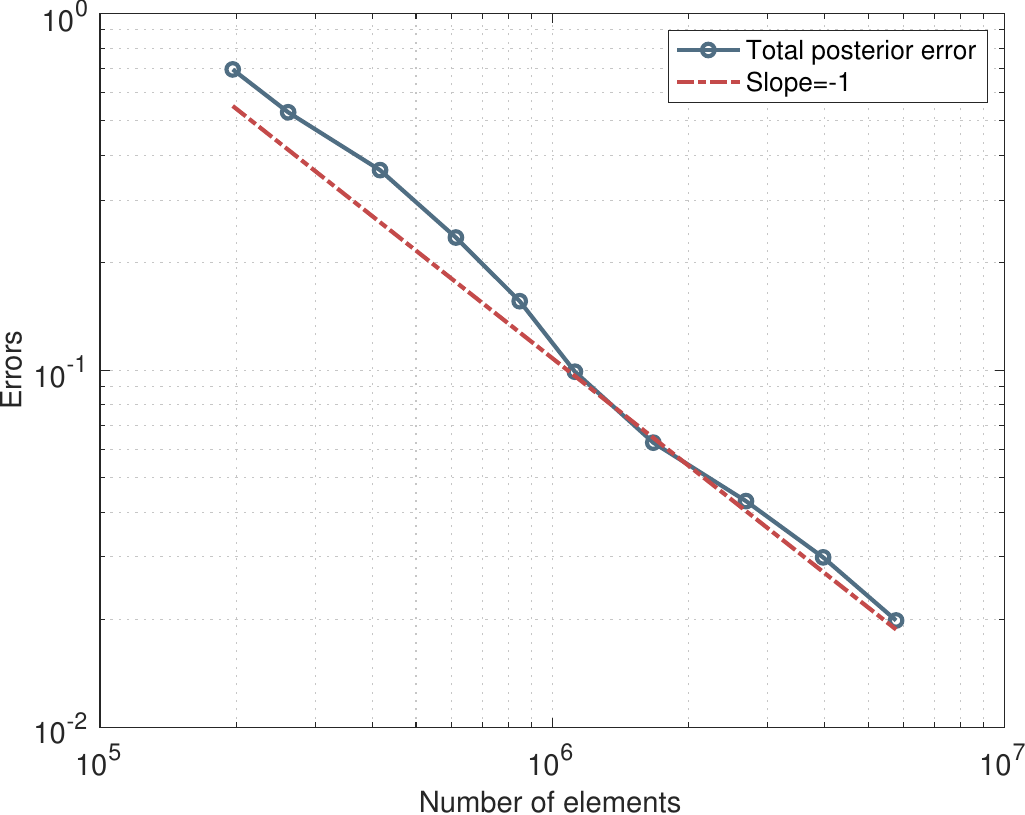}
  \caption{The convergence behavior of Algorithm \ref{Algorithm_PASE_Algebric} for the Hydrogen atom example, with $\rm nev = 200$ and $\rm tol = 10^{-8}$.}
  \label{fig:conv_ex4}
\end{minipage}
\hfill
\begin{minipage}{0.48\textwidth}
  \centering
  \includegraphics[width=\textwidth]{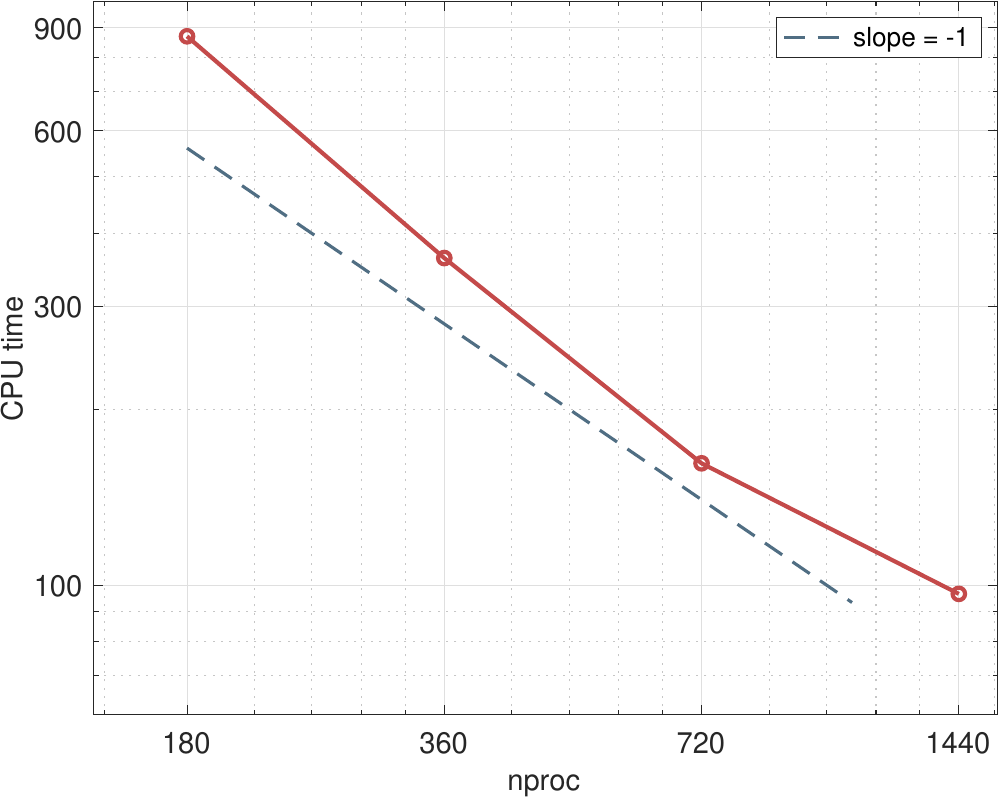}
  \caption{The scalability of Algorithm \ref{Algorithm_PASE_Algebric} for the Hydrogen atom example, with $\rm nev = 200$ and $\rm tol = 10^{-8}$.}
  \label{fig:scaling_ex4}
\end{minipage}
\end{figure}

The convergence behavior of Algorithm \ref{Algorithm_PASE_Algebric} for the Hydrogen atom example is shown in Figure \ref{fig:conv_ex4}. After 10 rounds of adaptive refinement, the number of elements reaches $N = 5,767,405$. The scalability of PASE is demonstrated in Figure \ref{fig:scaling_ex4}, where good scalability is observed for this eigenvalue problem.

\subsubsection{Example: L-Shaped Region}
In this example, we consider the eigenvalue problem defined on a non-convex region with singular properties:
Find $(\lambda,u)\in\mathcal R\times V$ such that $\|u\|_a=1$ and
\begin{eqnarray*}
-\Delta u=\lambda u, \quad {\rm in} \  \Omega,
\end{eqnarray*}
where $\Omega=(([-1,1]^2 / [0,1]^2))\times[0,1]$. It is important to note that the exact eigenpairs are obtained by solving the eigenvalue problem directly on a finely discretized finite element space.

The element residual for this L-shaped region example is defined as:
\begin{eqnarray*}\label{residual_Lshape}
\mathcal{R}_T\left(\Lambda_h, \Phi_h\right):=\left(\Delta u_{i, h}
-\lambda_{i, h} u_{i, h}\right)_{i=1}^N, 
\quad \text { in } T \in \mathcal{T}_{h}.
\end{eqnarray*}

We apply the adaptive mesh refinement strategy, similar to the Hydrogen atom example, to handle the singularities in the non-convex region. The jump residual part remains as defined earlier.

\begin{figure}[!htb]
\centering
\begin{minipage}{0.48\textwidth}
  \centering
  \includegraphics[width=\textwidth]{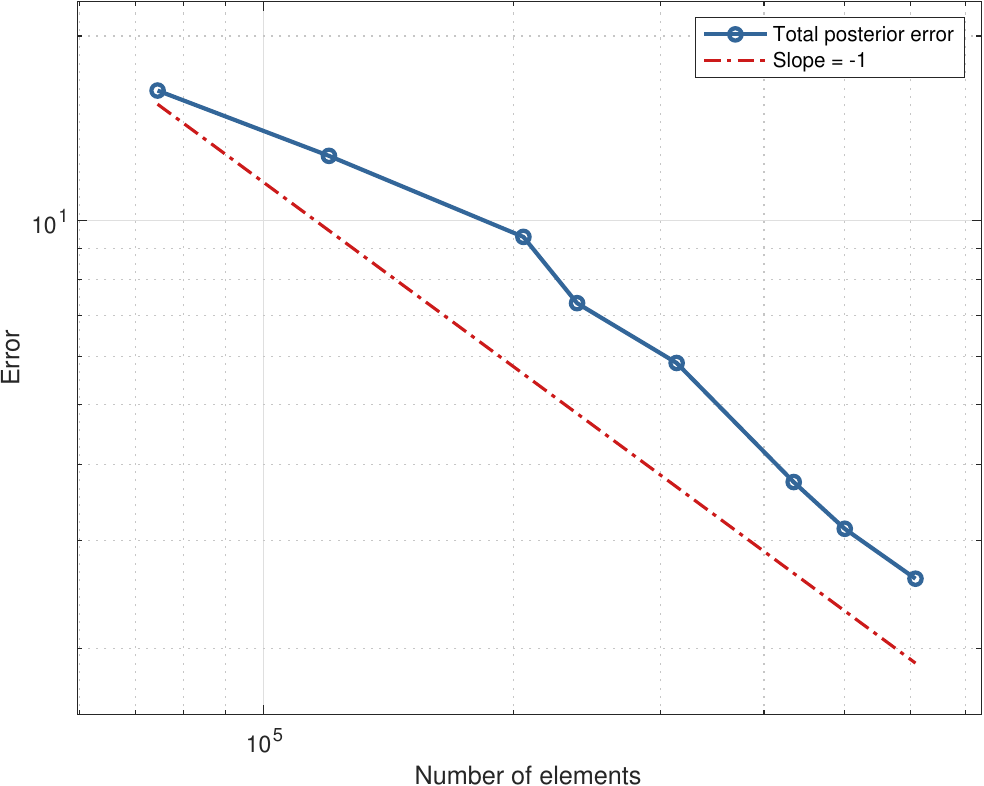}
  \caption{The convergence behavior of Algorithm \ref{Algorithm_k} for the L-Shaped Region example, with $\rm nev = 200$ and $\rm tol = 10^{-8}$.}
  \label{fig:Lshaped-ratio}
\end{minipage}
\hfill
\begin{minipage}{0.48\textwidth}
  \centering
  \includegraphics[width=\textwidth]{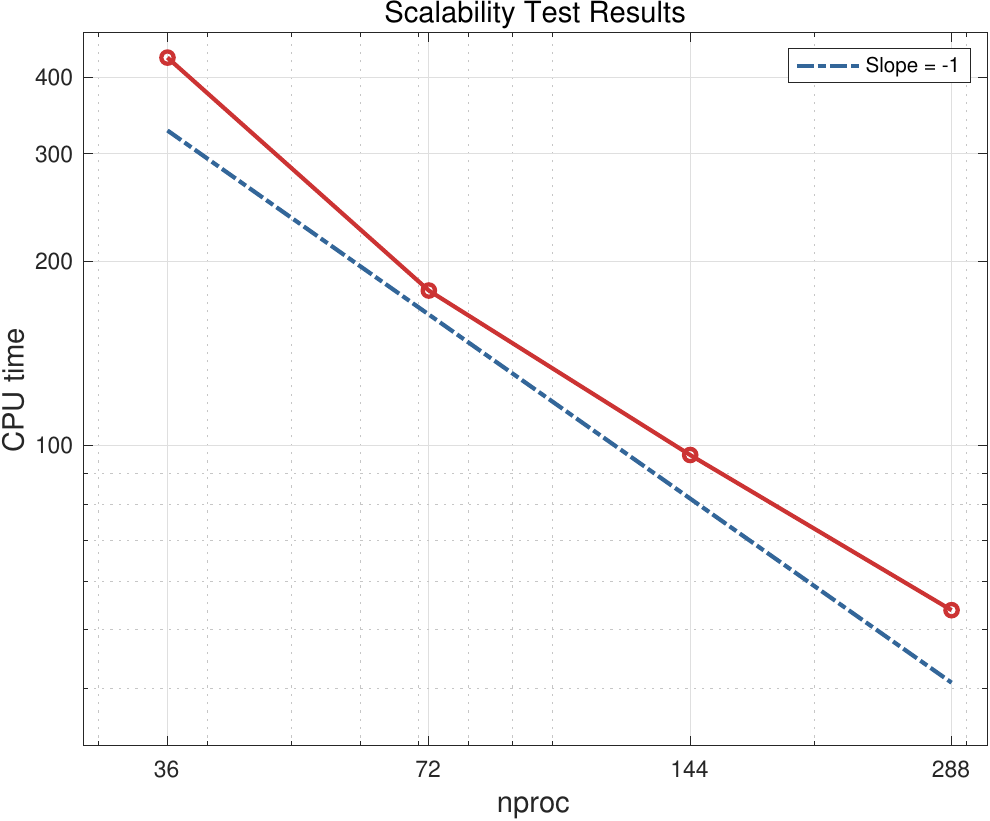}
  \caption{The scalability of Algorithm \ref{Algorithm_k} for the Hydrogen atom example, with $\rm nev = 200$ and $\rm tol = 10^{-8}$.}
  \label{fig:scaling_Lshaped}
\end{minipage}
\end{figure}
The convergence behavior of Algorithm \ref{Algorithm_PASE_Algebric} for the L-shaped region example is shown in Figure \ref{fig:Lshaped-ratio}. After 12 rounds of adaptive refinement, the number of elements reaches $N = 621,015$. The test uses $\rm nev = 200$ and $\rm tol = 10^{-8}$. The scalability of PASE is demonstrated in Figure \ref{fig:scaling_Lshaped}, where the algorithm continues to show excellent scalability despite the complexity introduced by the singularities in the region.

    \section{Conclusions}
    In this paper, we propose a type of parallel augmented subspace eigensover 
    for solving eigenvalue problems. In the parallel augmented subspace method, 
    the low dimensional augmented space is built with the finite element 
    space defined on the coarse mesh. 
    Enough numerical experiments are carried out to test the superiority of the 
    proposed parallel augmented subspace eigensolver to the GCGE, LOBPCG and Krylov-Schur 
    methods. Especially, compared with LOBPCG and Krylov-Schur, the proposed method has obviously better efficiency.



\end{document}